# Controllability on infinite-dimensional manifolds

M Khajeh Salehani and I Markina
Department of Mathematics, University of Bergen, NO-5008 Bergen, Norway

**Abstract.** Following the unified approach of A Kriegl and P W Michor (1997) for a treatment of global analysis on a class of locally convex spaces known as convenient, we give a generalization of Rashevsky–Chow's theorem for control systems in regular connected manifolds modelled on convenient (infinite-dimensional) locally convex spaces which are not necessarily normable.

## 1   Introduction

Control theory is in fact the theory of prescribing motion for dynamical systems rather than describing their observed behaviour.

The theory, originally developed to satisfy the design needs of servomechanisms, under the name of "automatic control theory," became recognized as a mathematical subject in 1956, with the publication of the paper of Boltyanski *et al* (1956) followed by the early papers of Kalman (1960, 1963) and Kalman *et al* (1962). Kalman challenged the accepted approach to control theory of that period, limited to the use of Laplace transforms and the frequency domain, by showing that the basic control problems could be studied effectively through the notion of a state of the system that evolves in time according to ordinary differential equations in which controls appear as parameters. Aside from drawing attention to the mathematical content of control problems, Kalman's work (1960) served as a catalyst for further growth of the subject. Liberated from the confines of the frequency domain and further inspired by the development of computers, automatic control theory became the subject matter of a new science called systems theory.

The initial orientation of systems theory, characteristic of the early 1960s, led away from geometric interpretations of linear theory and was partially responsible for the indifference with which Hermann's pioneering work (1963) relating Rashevsky–Chow's theorem to control theory was received by the mathematical community.

The significance of the Lie bracket for problems of control became clear around the year 1970 with publication of the papers of Brockett (1972), Hermes (1974) and Lobry (1970, 1974), followed by the papers of Brunovsky (1976), Elliott (1971), Krener (1974), Sussmann (1973), and others. Thanks to that collective effort, differential geometry entered into an exciting partnership with control theory. Control theory, on the other hand, through its distinctive concern for time-forward evolution of systems, led to its own theorems, making the birth of geometric control theory. For recent accounts on the geometric theory of control systems, we refer the reader to (Agrachev and Sachkov 2004) and (Bullo and Lewis 2005).

E−mail: salehani.k h.m@gmail.com (M.K. Salehani)

One of the fundamental problems in control theory is that of controllability, the question of whether one can drive the system from one point to another with a given class of controls. A classical result in geometric control theory of finite-dimensional (nonlinear) systems is Rashevsky–Chow's theorem that gives a sufficient condition for controllability on any connected manifold of finite dimension. In other words, the classical Rashevsky–Chow's theorem, which is in fact the first theorem in subriemannian geometry, gives a global connectivity property of a subriemannian manifold. The classical result was proved independently and almost simultaneously by Rashevsky (1938) and Chow (1939).

A similar result obtained by Carathéodory (1909) for analytic *distributions* of codimension one, in connection with his studies on the foundations of thermodynamics, has been extended in (Chow 1939) to smooth distributions of arbitrary codimension. Rashevsky (1938) was probably inspired by the vigorous research which was centred at that time in the seminars of Kagan, Finikov, and Vagner (1935). The classical theorem of Rashevsky and Chow was later proved by Sussmann (1973) under weaker conditions on the distributions -as compared to the *completely nonholonomic* condition. A distribution here means a linear subbundle of the tangent bundle of a manifold (Montgomery 2002).

It is worth noting here that there is a close link between nonholonomic constraints and controllability of nonlinear systems. Nonholonomic constraints are given by *nonintegrable* distributions, i.e., taking the Lie bracket of two vector fields in such a distribution may give rise to a vector field not contained in the same distribution. It is precisely this property which is needed in a nonlinear control system so that we can drive the system to as large a part of the state manifold as possible.

In our study of control systems, we always assume that the state space $M$ is a smooth manifold modelled on a locally convex space, and deal with the flows of some family $\mathcal{F} \subset \mathrm{Vec}(M)$ of complete smooth vector fields on $M$. Let $\mathcal{P}(\mathcal{F}) =: \mathcal{P}$ denote the group of diffeomorphisms of $M$ generated by flows $\{e^{tX} | t \in \mathfrak{R}\}_{X \in \mathcal{F}}$ of $\mathcal{F}$, $\mathrm{Lie}\,\mathcal{F}$ be the Lie subalgebra of $\mathrm{Vec}(M)$ generated by $\mathcal{F}$, and $\mathrm{Lie}_x \mathcal{F} = \{V(x) | V \in \mathrm{Lie}\,\mathcal{F}\}$ -the evaluation of $\mathrm{Lie}\,\mathcal{F}$ at $x \in M$. We say that $\mathcal{F} \subset \mathrm{Vec}(M)$ is *bracket generating*, or *completely nonholonomic*, if $\mathrm{Lie}_x \mathcal{F} = T_x M$, for every $x \in M$.

Accordingly, the classical version of Rashevsky–Chow's theorem states that: if $M$ is a connected manifold of finite dimension, and $\mathcal{F}$ is bracket generating then $\mathcal{P} \subset \mathrm{Diff}(M)$ acts transitively on $M$, i.e., $\mathcal{P}(x) = M$ for each $x \in M$.

This classical result, however, does not hold for infinite-dimensional control systems in general -i.e., the case when $M$ is of infinite dimension.

Some attempts have in fact been made to generalize the above-mentioned classical result to hold for infinite-dimensional state spaces. In their study of certain classes of "controllable" systems described by partial differential equations, Dudnikov and Samborski (1980) formulated a version of Rashevsky–Chow's theorem for control systems in any Banach vector spaces. In addition to the above-mentioned work, a generalization of Rashevsky–Chow's theorem for control systems in any complete connected *Hilbert manifold* (i.e., a manifold modelled on a Hilbert space) was given by Heintze and Liu (1999).

Now, the natural question arises whether it is possible to generalize Rashevsky–Chow's theorem for control systems in manifolds modelled on a more general class of topological vector spaces including the ones which, in contrast to Banach or Hilbert spaces, are not equipped with any norm or inner product.

In this article, following the convenient setting of infinite-dimensional differential geometry and global analysis developed by Kriegl and Michor (1997), we first introduce the notion of *Mackey completeness* in infinite-dimensional locally convex vector spaces in Section 2 which presents some preliminaries on a class of locally convex spaces known as *convenient*. Then in Section 3, we give a generalization of Rashevsky–Chow's theorem for control systems in *regular* connected manifolds $M$ modelled on convenient (infinite-dimensional) locally convex spaces $E$ (see Theorem 3.1); in fact, for a given family $\mathscr{F} \subset \mathrm{Vec}(M)$ of smooth vector fields on $M$, we show that if $\mathrm{Lie}_x \mathscr{F}$ is dense in $T_x M$ for all $x$ in $M$ then $\mathscr{P}(x)$ is dense in $M$ for all $x \in M$. We call here a smooth manifold regular if any neighbourhood of any $a \in M$ contains the closure of some smaller neighbourhood of the same point $a$ in $M$. The regularity condition is in fact satisfied if, for example, $M$ is locally compact or is a topological group (Milnor 1984).

In particular Theorem 3.1, which makes it possible to consider more general classes of controllable nonlinear systems including those of systems in Hilbert and Banach manifolds, gives also a refinement of Heintze–Liu's generalized version of Rashevsky–Chow's theorem because Heintze–Liu's conclusion of their theorem (1999) is not affected if we replace their completeness condition on Hilbert manifolds by Milnor's topological regularity condition introduced above.

The proof of Theorem 3.1 consists in the construction of some kind of *cones* in the locally convex vector spaces. The main difficulty in carrying out this construction is that locally convex vector spaces in general fail to have any norm or inner product; Lemma 2.12 is the key to constructing the cones. In fact, Lemma 2.12 and Corollary 3.2 are cornerstones of the proof of Theorem 3.1.

To indicate an application of Theorem 3.1, we conclude the paper with a controllability result on the group of orientation-preserving diffeomorphisms $\mathrm{Diff}_0(S^1)$; it is worth noting that this result does not follow from those obtained by Agrachev and Caponigro (2009), see Example 3.3.

## 2    Foundations: convenient locally convex spaces

In fact, classical calculus works quite well up to and including Banach spaces. There are many interesting works which have treated of global analysis mainly on manifolds modelled on Banach spaces; see for instance (Palais 1968) and (Eells 1966). However, further development has shown that Banach manifolds are not suitable for some classes of control systems and for infinite-dimensional global analysis, because as shown in (Omori and de la Harpe 1972) and (Omori 1978): if a Banach Lie group acts effectively on a finite-dimensional compact smooth manifold it must be finite dimensional itself. Moreover, Banach manifolds turn out to be open subsets of the modelling space in many cases, cf. (Eells and Elworthy 1970).

In his careful exposition of the Nash-Moser inverse function theorem, Hamilton (1982) defined a category of "tame" Fréchet spaces and investigated the manifolds modelled on them.

Differential calculus in infinite dimensions has already quite a long history; in fact, it goes back to Bernoulli and Euler, to the beginnings of variational calculus. During the twentieth century the urge to differentiate in spaces, which are more general than Banach and Fréchet spaces and are not necessarily normable, became stronger and many different approaches and definitions were attempted; e.g., a theory of differentiation was constructed by Yamamuro (1979) on locally convex spaces based on the correspondence between the sets of seminorms which induce original topologies.

To study locally convex spaces and the manifolds modelled on them, we follow the unified approach of Kriegl and Michor (1997) whose purpose is to lay the foundations of infinite-dimensional differential geometry on manifolds modelled on a class of locally convex spaces known as convenient.

We begin by introducing the required terminology using a sequence of definitions.

**Definition 2.1.** Let $E$ be a real vector space.
- A map $p: E \to \Re$ is said to be a *quasi-seminorm*, if
  (i) $p(x+y) \leq p(x) + p(y)$, for all $x, y \in E$;
  (ii) $p(tx) = tp(x)$, for all $x \in E$ and all $t \in \Re$ with $t \geq 0$.
- A map $p: E \to \Re$ is said to be a *seminorm* if, in addition to the above two properties, it satisfies:
  (ii') $p(\lambda x) = |\lambda| p(x)$, for all $x \in E$ and for all $\lambda \in \Re$.

It is evident that if $p: E \to \Re$ is a seminorm then $p(x) \geq 0$, for all $x \in E$.
(as $p(0) = p(x + (-x)) \leq p(x) + p(-x) = 2p(x)$ and $p(0) = 0$.)

The following proposition describes a method of constructing quasi-seminorms.

**Proposition 2.2.** Let $E$ be a real vector space. Suppose $D \subset E$ is a convex subset containing $0$, which is *absorbing*, i.e.,

$$\bigcup_{t>0} tD = E. \tag{1}$$

For every $x \in E$ we define

$$P_D(x) = \inf\{t > 0 \mid x \in tD\}. \tag{2}$$

(The set at the right-hand side of $(2)$ is nonempty since $D$ is absorbing.)

Then the map $P_D: E \to \Re$ is a quasi-seminorm; cf. (Jarchow 1981).

**Definition 2.3.** Under the hypothesis of the above proposition, the quasi-seminorm $P_D$ is called the *Minkowski functional* associated with the set $D$.

**Definition 2.4.** A *real topological vector space* is a vector space $E$, which is also a topological space, such that the maps

$$E \times E \ni (x, y) \mapsto x + y \in E$$
$$\Re \times E \ni (t, x) \mapsto tx \in E$$

are continuous.

**Lemma 2.5.** Let $E$ be a real topological vector space. Suppose $\tilde{D} \subset E$ is a convex open subset, which contains $0$. Then $\tilde{D}$ is absorbing, cf. $(1)$.

Moreover, the Minkowski functional associated with $D = -\tilde{D} \cap \tilde{D}$ is in fact a seminorm, where $-\mathcal{D} = \{-x | x \in \mathcal{D}\}$.

*Proof of Lemma 2.5.* To prove that $\bigcup_{t>0} t\mathcal{D} = E$, we define for each $x \in E$, the function
$$F_x : \mathfrak{R} \ni t \mapsto tx \in E.$$
Since $E$ is a topological vector space, the maps $F_x$ for all $x \in E$ are continuous. We start with an arbitrary $x \in E$, and use the continuity of the map $F_x$ at $t = 0$. Since $\mathcal{D}$ is a neighbourhood of $0$, there exists some $s > 0$ such that
$$F_x(t) \in \mathcal{D}, \ \forall t \in [-s, s].$$
In particular, $sx \in \mathcal{D}$ which means that $x \in s^{-1}\mathcal{D}$.

Let $D := -\mathcal{D} \cap \mathcal{D}$. It is clear that $-\mathcal{D}$ is a convex open set containing $0$, and so is the set $D$. Moreover, $D = -D$ from which it follows that the Minkowski functional $P_D$ is in fact a seminorm. Q.E.D.

***Definition 2.6.*** A topological vector space $E$ is said to be *locally convex* if for every $x \in E$ and every neighbourhood $U$ of $x$ there exists a convex open set $D$ such that $x \in D \subset U$.

***Definition 2.7.*** A set $M \subset E$ is called *bounded* if it is "absorbed" by each $0$-neighbourhood, i.e., for any neighbourhood $U \subset E$ of $0$ there exists a real number $\lambda > 0$ such that $M \subset \lambda U$.

***Definition 2.8.*** A set $M \subset E$ is *absolutely convex* if
$$\forall x_1, x_2 \in M, \ \{\lambda_1 x_1 + \lambda_2 x_2 \mid \lambda_i \in \mathfrak{R}, |\lambda_1| + |\lambda_2| \le 1\} \subset M.$$

For convenience of the reader we mention the following geometric version of the Hahn-Banach theorem without proof, thus making our exposition self-contained.

***Lemma 2.9.*** (Hahn-Banach Separation Theorem for Locally Convex Spaces). Let $E$ be a real locally convex vector space, and suppose $\mathcal{A}, \mathcal{B} \subset E$ are disjoint convex sets, with $\mathcal{A}$ compact, and $\mathcal{B}$ closed. Then there exists a linear continuous map $\ell : E \to \mathfrak{R}$, and two real numbers $\alpha, \beta \in \mathfrak{R}$, such that
$$\ell(x) \le \alpha < \beta \le \ell(y), \ \forall x \in \mathcal{A}, y \in \mathcal{B},$$
where $\alpha = \sup_{x \in \mathcal{A}} \ell(x)$, $\beta = \inf_{y \in \mathcal{B}} \ell(y)$.

We introduce the notion of *completeness* in *infinite-dimensional* locally convex vector spaces, following Kriegl and Michor (1997).

In finite-dimensional analysis we use the Cauchy condition, as a necessary condition for convergence of a sequence, to define completeness of our spaces. Here, in the infinite-dimensional case, we introduce the much stronger properties of being *Mackey-Cauchy* and being *Mackey convergent* as follows.

***Definition 2.10.*** Let $E$ be a locally convex space.
- A net $(x_\gamma)_{\gamma \in \Gamma}$ in $E$ is called *Mackey-Cauchy* provided that there exists a bounded absolutely convex set $M \subset E$ and a net $(\mu_{\gamma, \gamma'})_{(\gamma, \gamma') \in \Gamma \times \Gamma}$ in $\mathfrak{R}$ converging to $0$ such that $(x_\gamma - x_{\gamma'}) \in \mu_{\gamma, \gamma'} M$.

- For any bounded absolutely convex set $M \subset E$, we denote by $E_M$ the linear span of $M$ equipped with the Minkowski functional $P_M$, which is in fact a normed space. A net $(x_\gamma)_{\gamma \in \Gamma}$ is said to be convergent to $x$ in the normed space $(E_M, P_M)$ if there exists a net $\mu_\gamma \to 0$ in $\Re$ such that $x_\gamma \in \mu_\gamma M$.

- A net $(x_\gamma)_{\gamma \in \Gamma}$ in $E$ for which there exists a bounded absolutely convex $M \subset E$ such that $x_\gamma$ converges to $x$ in $E_M$ is called *Mackey convergent* (or briefly, *M-convergent*) to $x$.

- The space $E$ is said to be *Mackey complete* if every Mackey-Cauchy net converges in $E$.

Note that the above definition can also be given for sequences, in place of nets, in $E$ with a countable index set $\Gamma$.

The following result states when we call a vector space *convenient*.

**Lemma 2.11.** (Convenient Vector Spaces). Let $E$ be a locally convex vector space. $E$ is said to be *convenient* if one of the following equivalent (completeness) conditions is satisfied:
  (i)   $E$ is Mackey complete; i.e., every Mackey-Cauchy sequence converges.
  (ii)  If $M \subset E$ is absolutely convex closed bounded, then $E_M$ is a Banach space. This property is called *locally completeness* in (Jarchow 1981, p 196).
  (iii) For every bounded set $M \subset E$ there exists an absolutely convex bounded set $M' \supseteq M$ such that $E_{M'}$ is a Banach space.

The key to formulating the main results of this paper is the following lemma.

**Lemma 2.12.** Let $E$ be a convenient real locally convex vector space, and $B \subset E$ be a closed nonempty proper subset. Then there exists an element $a_* \in B$, an open set $U \subset E$ containing $a_*$, and a cone
$$C_{a_*} := \{a_* + t(x - a_*) \mid x \in \mathcal{X}, t \geq 0\},$$
such that $U \cap C_{a_*} \cap B = \{a_*\}$, where $\mathcal{X} \subset E$ is a convex closed set.

Before starting our proof of the above lemma, it is worth pointing out that for the case when $\dim(E) < \infty$, the above lemma can be rephrased by taking some ball (in place of the cone) with a point of $B$ on the boundary of the ball.

*Proof of Lemma 2.12.* Let $a_1 \in B$, $x_\circ \in E \setminus B$, and $D = -\mathcal{D} \cap \mathcal{D}$ with $\mathcal{D} \subset E$ being a proper convex open bounded set containing $0$. Lemma 2.5 now shows that $D$ is absorbing, and the Minkowski functional $P_D$ is a seminorm.

Since $E \setminus B \subset E$ is an open set in a locally convex vector space with $x_\circ \in E \setminus B$, there exists some convex open $V \subset E \setminus B$ with $x_\circ \in V$. Set $\alpha := \inf_{x \in \partial V} P_D(x - x_\circ)$, and define the

bounded set
$$S_{\alpha/2} := \{x \in V \mid P_D(x - x_\circ) \leq \alpha/2\}. \tag{3}$$

It is worth noting that $\alpha > 0$ since $x_\circ \notin \partial V$, and that $\alpha < \infty$ because $D$ is absorbing and $\partial V \neq \emptyset$. It is in fact immaterial which $D$ we choose to define $\alpha$ as long as $D$ is absorbing, and therefore $P_D$ always returns finite numbers.

Since $P_D$ is a seminorm, from (2) and Definitions 2.1 and 2.10, it follows that

**Fact 2.12.1.** $S_{\alpha/2}$ is convex, and so is its closure $\bar{S}_{\alpha/2}$ in $E$; and that

**Fact 2.12.2.** $\bar{S}_{\alpha/2}$ and $\partial V$ are disjoint sets, and consequently $\bar{S}_{\alpha/2} \cap B = \emptyset$.

Now, define a solid cone with vertex at $a_1 \in B$ as
$$C_{a_1} := \{a_1 + t(x - a_1) \mid x \in \bar{S}_{\alpha/4}, t \geq 0\},$$
where $S_{\alpha/4} = \{x \in V \mid P_D(x - x_\circ) \leq \alpha/4\} \subset S_{\alpha/2} \subset V$. It follows from Fact 2.12.2 that $\bar{S}_{\alpha/4} \cap B = \emptyset$. Consider
$$B_1 := \{a_1 + t(x - a_1) \mid x \in \bar{S}_{\alpha/4}, 0 \leq t \leq 1\} \cap B.$$
It is immediate that $B_1 \subset C_{a_1} \cap B$ is a closed bounded set in $E$.

In fact, the proof of Lemma 2.12 is based on the following claim.

**Claim 2.12.3.** There exists an element $a_* \in B$ such that
$$C_{a_*} \cap B_1 = \{a_*\},$$
where $C_{a_*}$ is a cone with vertex at $a_*$.

*Proof of Claim 2.12.3.* The Hahn-Banach separation theorem for locally convex vector spaces (cf. Lemma 2.9), with $\mathscr{A} := \{a_1\}$ and $\mathscr{B} := \bar{S}_{\alpha/4}$, shows that there exist a linear continuous map $\ell : E \to \Re$ and a real number $\beta \in \Re$ such that
$$\ell(a_1) < \beta \leq \ell(y), \quad \forall y \in \bar{S}_{\alpha/4}.$$
Let $\ell(x) = c$ (for some real number $\ell(a_1) < c < \beta$) be a hyperplane in $E$ separating $a_1 \in B$ and $\bar{S}_{\alpha/4}$. Set $e := \dfrac{x_\circ - a_1}{\ell(x_\circ - a_1)}$, and consider the set of points $\{a_1 + (\ell(b) - \ell(a_1))e \mid b \in B_1\}$ which can be thought of as the "projection" of $B_1$ on the $e$-axis (i.e., on the 1-dimensional affine subspace $a_1 + \Re e$ in $E$).

It is evident that $\ell(b) - \ell(a_1) \geq 0$, for all $b \in B_1$, so $B_1$ is projected on the positive half of the $e$-axis and
$$d := \sup_{b \in B_1} (\ell(b) - \ell(a_1))$$
is non-negative. If $d = 0$ then $B_1$ is projected to $\{a_1\}$ and hence $C_{a_1} \cap B_1 = \{a_1\}$, so we are done.

Now suppose that $d > 0$, and set
$$C_{a_1}^{\Pi_{a_1}} := C_{a_1} \cap \{x \in E \mid \ell(x) \leq \ell(a_1) + d\},$$

which is the cone $C_{a_1}$ truncated by the hyperplane $\Pi_{a_1} : \ell(x) = \ell(a_1) + d$.

Clearly the closed set $C_{a_1}^{\Pi_{a_1}} \subset \{a_1 + t(x - a_1) | x \in \bar{S}_{\alpha/2}, 0 \leq t \leq 1\}$ is bounded. Since $d > 0$, it follows that $C_{a_1}^{\Pi_{a_1}} \cap B_1 \supset \{a_1\}$ and hence there exists $a_2 \in C_{a_1}^{\Pi_{a_1}} \cap B_1$ such that
$$\ell(a_2) - \ell(a_1) > d/2. \tag{4}$$

Denote the parallel translation of the cone $C_{a_1}$ along the vector $(a_2 - a_1)$ by $C_{a_2}$ and define $B_2 := C_{a_2} \cap B_1$.

If $B_2 = C_{a_2} \cap B_1 = \{a_2\}$ then we are done. Now suppose that $B_2 \supset \{a_2\}$ and define the truncated cone
$$C_{a_2}^{\Pi_{a_1}} := C_{a_2} \cap \{x \in E | \ell(x) \leq \ell(a_1) + d\},$$
the latter being contained in $C_{a_1}^{\Pi_{a_1}}$.

Since $E$ is a convenient vector space and $C_{a_1}^{\Pi_{a_1}}$ is bounded, it follows from Lemma 2.11 that there exists an absolutely convex bounded set $C' \supseteq C_{a_1}^{\Pi_{a_1}}$ such that $(E_{C'}, P_{C'})$ is a Banach space; indeed, $C'$ can be equal to $\lambda D$ for some real number $\lambda > 0$ because $C_{a_1}^{\Pi_{a_1}}$ is bounded and $D$ is a $0$-neighbourhood, cf. Definition 2.7.

The *diameter* of $C_{a_1}^{\Pi_{a_1}}$, denoted by $\mathrm{diam}(C_{a_1}^{\Pi_{a_1}})$, is defined to be
$$\mathrm{diam}(C_{a_1}^{\Pi_{a_1}}) := \sup_{x,y \in C_{a_1}^{\Pi_{a_1}}} P_{C'}(x - y). \tag{5}$$

The diameter of $C_{a_2}^{\Pi_{a_1}}$ and that of any other subset of $C' \supseteq C_{a_1}^{\Pi_{a_1}}$ can be defined similarly to above.

In order to compare the diameters of $C_{a_i}^{\Pi_{a_1}}$, for $i = 1, 2$, it is convenient to parallel translate the whole space $E$ along the vector $(a_1 - a_2)$, which gives another copy of $E$. Thus $C_{a_2}$ is transformed to $C_{a_1}$, and the set $\{x \in E | \ell(x) \leq \ell(a_1) + d\}$ coincides with $\{x \in E | \ell(x) \leq 2\ell(a_1) - \ell(a_2) + d\}$. Therefore $C_{a_2}^{\Pi_{a_1}}$ will be transformed to
$$C_{a_1}^{\Pi_{(2a_1 - a_2)}} := C_{a_1} \cap \{x \in E | \ell(x) \leq 2\ell(a_1) - \ell(a_2) + d\},$$
which is the cone $C_{a_1}$ truncated by the hyperplane $\Pi_{(2a_1 - a_2)} : \ell(x) = 2\ell(a_1) - \ell(a_2) + d$, and is in fact contained in $C_{a_1}^{\Pi_{a_1}}$ since $2\ell(a_1) - \ell(a_2) + d < \ell(a_1) + d$.

Let $(a_1 + t_1 e)$ be the intersection point of the hyperplane $\Pi_{a_1}$ with the 1-dimensional

(affine) subspace $(a_1+\Re e)\subset E$, and $(a_1+t_2 e)$ be that of the hyperplane $\Pi_{(2a_1-a_2)}$ with $(a_1+\Re e)$. Evidently, the sets $C_{a_1}^{\Pi_{(2a_1-a_2)}}$ and $C_{a_1}^{\Pi_{a_1}}$ are homothetic with the coefficient equals $k_1 = t_2/t_1$, where $t_1 = d$ and $t_2 = d-(\ell(a_2)-\ell(a_1))$. Therefore,

$$\rho_2 := \mathrm{diam}(C_{a_2}^{\Pi_{a_1}}) = \mathrm{diam}(C_{a_1}^{\Pi_{(2a_1-a_2)}}) = k_1 \mathrm{diam}(C_{a_1}^{\Pi_{a_1}}) =: k_1 \rho_1.$$

It follows from (4) that $\rho_2 = k_1 \rho_1 < \rho_1/2$.

Now, project $B_2 = C_{a_2} \cap B_1$ on $(a_1+\Re e)$ as above. If $C_{a_2}^{\Pi_{a_1}} \cap B_2 \supset \{a_2\}$ then we may choose $a_3 \in C_{a_2}^{\Pi_{a_1}} \cap B_2$ such that $\ell(a_3)-\ell(a_2) > d/2$, and define the truncated cone

$$C_{a_3}^{\Pi_{a_1}} := C_{a_3} \cap \{x \in E \mid \ell(x) \le \ell(a_1)+d\},$$

where $C_{a_3}$ is obtained by the parallel translation of $C_{a_2}$ along $(a_3-a_2)$ and, in this way, can also be defined as that of $C_{a_1}$ along $(a_3-a_1)$.

Similarly to above, it can be seen that

$$\rho_3 := \mathrm{diam}(C_{a_3}^{\Pi_{a_1}}) < \rho_2/2 < \frac{\rho_1}{2^2}.$$

By the above procedure, we can construct two sequences of nested closed bounded sets

$$B \supset B_1 \supset B_2 \supset \cdots$$
$$C_{a_1}^{\Pi_{a_1}} \supset C_{a_2}^{\Pi_{a_1}} \supset \cdots$$

where $B_n = C_{a_n} \cap B_1$, the cone $C_{a_n}$ being obtained by parallel translation of $C_{a_1}$ along $(a_n-a_1)$, the truncated cone $C_{a_n}^{\Pi_{a_1}} = C_{a_n} \cap \{x \in E \mid \ell(x) \le \ell(a_1)+d\}$ with $\mathrm{diam}(C_{a_{n+1}}^{\Pi_{a_1}}) < \frac{\rho_1}{2^n}$, and $\{a_n\}_1^\infty \subset B_1$ is a sequence with $a_{n+1} \in C_{a_n}^{\Pi_{a_1}} \cap B_n = C_{a_n}^{\Pi_{a_1}} \cap B_1$.

Thus $\{a_n\}_1^\infty \subset B_1$ is a Mackey-Cauchy sequence because for any given $r < n$,

$$a_{n+2-r} \in C_{a_{n+1-r}}^{\Pi_{a_1}} \cap B_1 \subseteq \frac{1}{2^{n-r}} C_{a_1}^{\Pi_{a_1}},$$

$$a_{n+2} \in C_{a_{n+1}}^{\Pi_{a_1}} \cap B_1 \subset C_{a_{n+1-r}}^{\Pi_{a_1}} \cap B_1 \subseteq \frac{1}{2^{n-r}} C_{a_1}^{\Pi_{a_1}}.$$

Hence

$$(a_{n+2}-a_{n+2-r}) \in \frac{2^r}{2^n} C'$$

where $C'$ is an absolutely convex bounded set, and $\frac{2^r}{2^n}$ converges to $0 \in \Re$.

Therefore $\{a_n\}_1^\infty \subset B_1$ is (Mackey-)convergent to some element $a_* \in B_1$. Denote the parallel translation of the cone $C_{a_1}$ along the vector $(a_*-a_1)$ by $C_{a_*}$. It is obvious that $a_* \in C_{a_*} \cap B_1 \subset C_{a_1}^{\Pi_{a_1}}$. Moreover, since $(E_{C'}, P_{C'})$ is a Banach space, it follows from
$$C_{a_*} \cap B_1 \subset C_{a_{n+1}} \cap B_1 \subset C_{a_{n+1}}^{\Pi_{a_1}}$$
that, for any natural number $n$,
$$\operatorname{diam}(C_{a_*} \cap B_1) \leq \operatorname{diam}(C_{a_{n+1}}^{\Pi_{a_1}}) < \frac{\rho_1}{2^n} \to 0$$
which finishes the proof of Claim 2.12.3. Q.E.D.

Now consider
$$S_{\alpha/3} := \{x \in V \mid P_D(x - x_\circ) \leq \alpha/3\},$$
the positive real number $\alpha$ being as in (3); then, on account of Fact 2.12.1, its closure $\overline{S}_{\alpha/3}$ and interior $\operatorname{int}(\overline{S}_{\alpha/3})$ are both convex. Since $\operatorname{int}(\overline{S}_{\alpha/3})$ is a proper subset of $\overline{S}_{\alpha/2}$, there exists some real $\varepsilon > 0$ such that
$$(\operatorname{int}(\overline{S}_{\alpha/3}) + \varepsilon e) \subset \overline{S}_{\alpha/2}, \tag{6}$$
where $e = \dfrac{x_\circ - a_1}{\ell(x_\circ - a_1)}$ as above.

Define $U := \{a_1 + t(x - a_1) \mid x \in (\operatorname{int}(\overline{S}_{\alpha/3}) + \varepsilon e), 0 < t < 1\}$. It is easily seen that $U \subset E$ is an open set which contains $a_* \in B_1$, and that $U \cap C_{a_1} \cap B = B_1$, the latter being due to Fact 2.12.2 and (6); see Figure 1.

Note that $B$ may consist of several components, which is a reason for taking the neighbourhood $U$ of $a_*$ into account.

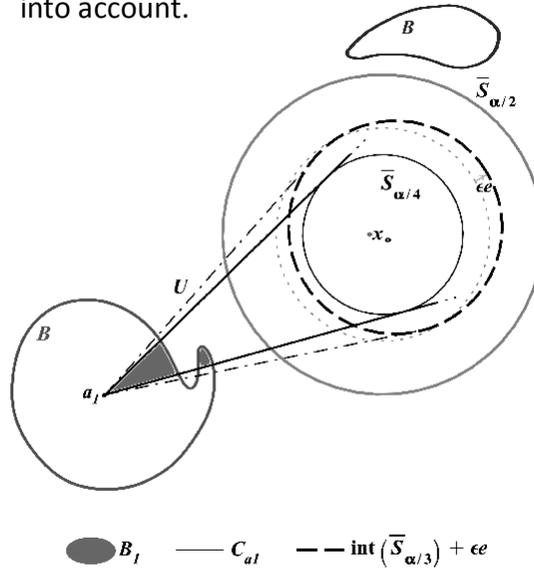

Figure 1. Illustrating the proof of Lemma 2.12.

Furthermore, since $U \cap C_{a_*} \cap B \subset U \cap C_{a_1} \cap B$, it follows that
$$U \cap C_{a_*} \cap B = (U \cap C_{a_1} \cap B) \cap C_{a_*} = B_1 \cap C_{a_*} = \{a_*\},$$
which completes the proof of Lemma 2.12. Q.E.D.

E−mail: salehani.k h.m@gmail.com (M.K. Salehani)

# 3   Controllability

One of the fundamental problems in control theory is that of controllability. Indeed, many design methodologies rely on some hypotheses that concern controllability. The problem of controllability is essentially that of describing the nature of the set of states reachable from an initial state. In the development of control theory, there are in fact two properties that arise as being important; namely, the property of "accessibility" and that of "controllability".

The property of accessibility means that the reachable set has a nonempty interior. The treatment of accessibility, in which we are interested, follows the approach of the fundamental paper of Sussmann and Jurdjevic (1972). Results of a related nature can be found in (Krener 1974) and (Hermann and Krener 1977). The property of controllability extends accessibility by further asking whether the initial state lies in the interior of the reachable set, i.e., the question of whether one can drive the system from one point to another with a given class of controls.

The matter of providing general conditions for determining controllability is currently unresolved, although there have been many deep and insightful contributions. Sussmann has made various important contributions to controllability, starting with the paper (Sussmann 1978). In the paper (Sussmann 1983), a Lie series approach was developed for the controllability of control-affine systems, and this approach culminated in the quite general results of (Sussmann 1987), which incorporated the ideas of Crouch and Byrnes (1986) concerning input symmetries. The Lie series methods rely on the notion that a system can be well approximated by a "nilpotent approximation." Contributions to this sort of approach have been made, for example, by Hermes (1982a,b, 1991), Kawski (1988b, 1998), and by Kawski and Sussmann (1997). A recent paper by Bianchini and Kawski (2003) indicates that there may well be some limitations to the approach of using nilpotent approximations to determine conditions for controllability. A related approach is the construction of "control variations," which is explained, for example, in the papers (Bianchini and Stefani 1993) and (Kawski 1988a).

Another approach to local controllability is that taken by Agrachev and Gamkrelidze (1993a,b), based on the chronological calculus of the same authors (Agrachev and Gamkrelidze 1978) and (Agrachev and Sachkov 2004). The fact that some of the very basic properties of the reachable set for a nonlinear control system are yet to be understood is the subject of the open problems paper by Agrachev (1999). Sontag (1988) and Kawski (1990) showed that a general answer to the controllability problem will be computationally difficult. Nonetheless, the problem of controllability is so fundamental that there continues to be much work in the area.

Other treatments of nonlinear controllability, in textbook form, include (Agrachev and Sachkov 2004), (Isidori 1995), (Jurdjevic 1997), (Nijmeijer and van der Schaft 1990) for accessibility, and (Bloch 2003) for accessibility and controllability. Some global controllability results are given in (Agrachev and Sachkov 2004), (Agrachev and Caponigro 2009) and (Grong *et al* 2012a,b).

The approach we follow here is based on the ones initiated in the works of Jurdjevic, Agrachev and Sachkov; for a through treatment, we refer the reader to (Jurdjevic 1997) and (Agrachev and Sachkov 2004) and the references given there.

In order to study infinite-dimensional smooth manifolds $M$ modelled on convenient locally convex spaces $E$, we need to give a brief exposition of the notion of smoothness for mappings between such manifolds and introduce the *kinematic* tangent bundles and vector fields; for a complete account of the infinite-dimensional differential geometry on such manifolds, we refer the reader to (Kriegl and Michor 1997).

Since the notion of smooth curves can be given without problems, a mapping between smooth manifolds modelled on convenient locally convex spaces is said to be smooth if it maps smooth curves into smooth curves. This notion of smoothness coincides with the usual concepts, up to manifolds modelled on Fréchet spaces.

For any $x \in M$, we say that $v$ is a *kinematic* tangent vector to $M$ at $x$ if there exists a curve $\gamma_v : [0,1] \to M$ such that $\gamma_v(0) = x$ and $\dot{\gamma}_v(0)$ exists and is equal to $v$, which explains the choice of the name kinematic. The kinematic tangent space of $M$ at $x$, denoted by $T_x M$, is defined to be the set of all kinematic tangent vectors to $M$ at $x$. As is the case for manifolds of finite dimensions, a chart map induces a one-to-one correspondence between the model space $E$ and a kinematic tangent space of $M$. Using these one-to-one correspondences, the kinematics tangent spaces can evidently be given the structure of topological vector spaces isomorphic to the convenient locally convex space $E$. Similarly the disjoint union of the kinematic tangent spaces $T_x M$, as $x$ varies over $M$, can be made into a new smooth manifold $TM$, which is called the kinematic tangent bundle of $M$ and is modelled on the locally convex space $E \times E$.

A kinematic vector field on $M$ is just a smooth section of the kinematic tangent bundle $TM$. In fact, some of the classically equivalent definitions of tangent vectors differ in infinite dimensions, and accordingly we have two different kinds of tangent bundles and vector fields; namely the "operational" ones, and those of the kinematic type introduced above. However, throughout this paper, we will be concerned only with the kinematic type because only kinematic vector fields can have flows which are in fact unique if they exist.

The control systems that we consider here will always be of the following form.

The state space $M$ is a smooth manifold modelled on a locally convex space, the control set $\mathcal{U}$ is an arbitrary (usually closed) subset of some Euclidean space, and the dynamics are described by a mapping $F : M \times \mathcal{U} \to TM$ such that for each $u \in \mathcal{U}$, $F_u : M \to TM$ defined by $F_u(x) = F(x,u) \in T_x M$ for $x$ in $M$ is a smooth vector field. Setting $\mathcal{F} := \{F_u | u \in \mathcal{U}\}$ to be the family of vector fields generated by $F$, we call a continuous curve $x : [0,T] \to M$ an *integral curve of* $\mathcal{F}$ if there exist a partition $0 = t_0 < t_1 < \cdots < t_m = T$ and vector fields $X_1, \ldots, X_m$ in $\mathcal{F}$ such that the restriction of $x(t)$ to each open interval $(t_{i-1}, t_i)$ is differentiable, and $dx(t)/dt = X_i(x(t))$ for $i = 1, \ldots, m$. In fact, $x(t)$ can be visualized as a "broken" continuous curve consisting of pieces of integral curves of vector fields corresponding to different choices of control values.

In what follows, $\mathcal{F} \subset \text{Vec}(M)$ stands for any family of smooth vector fields. To simplify the notations, we assume that all vector fields in $\mathcal{F}$ are complete. Thus each element $X \in \mathcal{F}$ generates a one-parameter group of diffeomorphisms $\{e^{tX} | t \in \Re\}$ = flow of $X$ in $M$. Let $\mathcal{P}(\mathcal{F}) =: \mathcal{P}$ denote the group of diffeomorphisms of $M$ generated by flows $\{e^{tX} | t \in \Re\}_{X \in \mathcal{F}}$ of $\mathcal{F}$. Each element $\Phi$ of $\mathcal{P} \subset \text{Diff}(M)$ is of the form
$$\Phi = e^{t_k X_k} \circ e^{t_{k-1} X_{k-1}} \circ \cdots \circ e^{t_1 X_1},$$

for some natural number $k$, real numbers $t_1,\ldots,t_k \in \Re$ and some vector fields $X_1,\ldots,X_k \in \mathcal{F}$. $\mathcal{P}(\mathcal{F}) = \mathcal{P}$ acts on $M$ in the obvious way and partitions $M$ into the sets $\mathcal{P}(x) = \{\Phi(x) | \Phi \in \mathcal{P}\}$ for $x$ in $M$.

Since the set $\text{Vec}(M)$ of all smooth vector fields on $M$ has the structure of a real Lie algebra under the Lie-bracket operation, to the given $\mathcal{F} \subset \text{Vec}(M)$ there corresponds the Lie subalgebra $\text{Lie}\,\mathcal{F}$ of $\text{Vec}(M)$ generated by $\mathcal{F}$. The evaluation of $\text{Lie}\,\mathcal{F}$ at $x \in M$ will be denoted by $\text{Lie}_x\,\mathcal{F} = \{V(x) | V \in \text{Lie}\,\mathcal{F}\}$.

A family $\mathcal{F} \subset \text{Vec}(M)$ is called *bracket generating*, or *completely nonholonomic*, if

$$\text{Lie}_x\,\mathcal{F} = T_xM, \text{ for every } x \in M. \tag{7}$$

A classical result in geometric control theory of finite-dimensional nonlinear systems, which gives a sufficient condition for controllability, is Rashevsky–Chow's theorem; cf. (Rashevsky 1938) and (Chow 1939).

If $M$ is a connected manifold of *finite* dimension, and $\mathcal{F}$ is bracket generating then Rashevsky–Chow's theorem states that the group $\mathcal{P}(\mathcal{F}) = \mathcal{P} \subset \text{Diff}(M)$ acts transitively on $M$, i.e., $\mathcal{P}(x) = M$ for each $x$ in $M$.

In fact, the relevance of the Lie bracket and Frobenius' theorem (1911) for controllability studies of finite-dimensional nonlinear systems of the form $dx/dt = F(x,u)$ comes in via the theorem of Rashevsky and Chow, and its refinement by others (Hermann 1963), (Hermes and Haynes 1963), (Krener 1974) and (Sussmann and Jurdjevic 1972).

The following theorem is a generalization of the above classical result to the case of *infinite*-dimensional manifolds, which makes it possible to consider even more general classes of "controllable" nonlinear systems.

To formulate the following theorem we need to introduce the notion of *topological regularity* of smooth manifolds. We call a smooth manifold *regular* if for any neighbourhood $V$ of any $a \in M$ there exists a neighbourhood $V' \subset M$ of $a$ such that its closure $\overline{V'} \subset V$. The regularity condition is in fact satisfied if, for example, $M$ is locally compact or is a topological group (Milnor 1984, p 1029).

**Theorem 3.1.** Let $M$ be a regular connected manifold modelled on a convenient locally convex space $E$, and $\mathcal{F}$ be a family of smooth vector fields on $M$. If $\text{Lie}_x\,\mathcal{F}$ is dense in $T_xM$ for all $x$ in $M$, then $\mathcal{P}(x)$ is dense in $M$ for all $x \in M$.

The remainder of this section will be devoted to a proof of the above theorem.

Let $\mathcal{B}$ be an arbitrary subset of the manifold $M$. For any $x \in \mathcal{B}$ and $v \in T_xM$, we say that $v$ is *tangent to* $\mathcal{B}$ *at* $x$ if there exists a curve $\gamma_v : [0,1] \to M$ such that $\gamma_v(0) = x$, $\dot{\gamma}_v(0)$ exists and is equal to $v$, and $\gamma_v(t) \in \mathcal{B}$ for all $t$. We denote by $\mathcal{T}_x\mathcal{B}$ the set of all tangent vectors to $\mathcal{B}$ at $x$.

The proof of Theorem 3.1 is based on the following corollary of Lemma 2.12.

***Corollary 3.2.*** Let $M$ be a regular connected manifold modelled on a convenient locally convex space $E$, and $\mathcal{B} \subseteq M$ be a closed nonempty subset. If $\mathcal{T}_x \mathcal{B}$ is dense in $T_x M$ for every $x \in \mathcal{B}$ then $\mathcal{B} = M$.

*Proof of Corollary 3.2.* Since $M$ is connected, we only need to prove that $\mathcal{B}$ is also open. On the contrary, suppose that $\mathcal{B}$ is *not* open. Then there exists a boundary point $b$ of $\mathcal{B}$. Let $(\varphi, V)$ be a chart around $b$. It follows from the regularity of the smooth manifold $M$ that there exists a neighbourhood $V' \subset M$ of $b$ such that $\overline{V'} \subset V$. Hence $\varphi(\mathcal{B} \cap \overline{V'}) \subset E$ is a closed nonempty proper subset.

Lemma 2.12 now shows that there exists an element $p \in \varphi(\mathcal{B} \cap \overline{V'})$, an open set $U \subset E$ containing $p$, and a cone $C_p$ with vertex at $p$ such that $U \cap C_p \cap \varphi(\mathcal{B} \cap \overline{V'}) = \{p\}$. As in the proof of Lemma 2.12, $p$ is a boundary point of $\varphi(\mathcal{B} \cap \overline{V'})$. So there exists a sequence $\{p_n\}_1^\infty \subset \varphi(\mathcal{B} \cap \overline{V'})$ that is (Mackey-)convergent to $p$. Since the chart map $\varphi$ is a homeomorphism, the sequence $\{\varphi^{-1}(p_n)\}_1^\infty \subset \mathcal{B} \cap \overline{V'}$ converges to $\varphi^{-1}(p) \in \partial(\mathcal{B} \cap \overline{V'})$. On the other hand, $\mathcal{C}_{\varphi^{-1}(p)} \cap \varphi^{-1}(U) \cap (\mathcal{B} \cap \overline{V'}) = \{\varphi^{-1}(p)\}$ where $\mathcal{C}_{\varphi^{-1}(p)} := \varphi^{-1}(C_p)$ is the subset of $M$ diffeomorphic to the cone $C_p \subset E$.

Consequently,
$$\mathcal{T}_{\varphi^{-1}(p)} \mathcal{B} \subset \overline{\mathsf{T}_{\varphi^{-1}(p)} \mathcal{B}} \subset \overline{\mathcal{T}_{\varphi^{-1}(p)}(\mathcal{C}_{\varphi^{-1}(p)} \cup \mathcal{B})} \subseteq T_{\varphi^{-1}(p)} M,$$
which contradicts the assumption that $\mathcal{T}_x \mathcal{B}$ is dense in $T_x M$ for all $x \in M$.  Q.E.D.

It is worth pointing out that for any family of smooth vector fields $\mathcal{F} \subset \mathrm{Vec}(M)$,
$$\mathrm{Lie}_x \mathcal{F} \subseteq \mathcal{T}_x \mathcal{P}(x), \text{ for every } x \in M. \tag{8}$$
This is proved by taking the following steps.

If $X \in \mathcal{F}$ then $t \mapsto e^{tX}(x)$ is a trajectory in $\mathcal{P}(x)$ whose velocity vector $X(x)$, at $t = 0$, is in $\mathcal{T}_x \mathcal{P}(x)$. If we take two arbitrary vector fields $X, Y \in \mathcal{F}$ then the diffeomorphism $e^{-tY} \circ e^{-tX} \circ e^{tY} \circ e^{tX}$ is in $\mathcal{P}$. In fact, the vector $[X,Y](x)$ is tangent to the trajectory $t \mapsto (e^{-tY} \circ e^{-tX} \circ e^{tY} \circ e^{tX})(x) \in \mathcal{P}(x)$ at $t = 0$, i.e., $[X,Y](x) \in \mathcal{T}_x \mathcal{P}(x)$. The rest of the proof runs, as above, by induction on the natural number $k$ as in the definition of $\mathrm{Lie}\,\mathcal{F} = \mathrm{span}\{[X_1, [\ldots [X_{k-1}, X_k] \ldots]] \mid X_i \in \mathcal{F}, k = 1, 2, 3, \cdots\} \subset \mathrm{Vec}(M)$.

*Proof of Theorem 3.1.* We first claim that
$$\mathcal{P}(x) \subset \overline{\mathcal{P}(p)}, \text{ for any } p \in M, x \in \overline{\mathcal{P}(p)}. \tag{9}$$
To see this, let $p \in M$ and $x \in \overline{\mathcal{P}(p)}$. Then there exists a sequence $\{x_n\}_1^\infty \subset \mathcal{P}(p)$ with $\lim_{n \to \infty} x_n = x$. For every $y \in \mathcal{P}(x)$ there exists a diffeomorphism $\Phi \in \mathcal{P} \subset \mathrm{Diff}(M)$ such that $y = \Phi(x)$. Moreover, to the sequence $\{x_n\}_1^\infty \subset \mathcal{P}(p)$ there corresponds a sequence of diffeomorphisms $\{\Psi_n\}_1^\infty \subset \mathcal{P}$ such that $x_n = \Psi_n(p)$. It follows that $\{\Phi(x_n)\}_1^\infty \subset \mathcal{P}(p)$, and $\lim_{n \to \infty} \Phi(x_n) = \Phi(x) = y \in \overline{\mathcal{P}(p)}$, and so (9) is proved.

By (8) and (9), $\text{Lie}_x \mathcal{F} \subseteq \mathcal{T}_x \mathcal{P}(x) \subset \mathcal{T}_x \overline{\mathcal{P}(p)} \subset T_x M$ for any $p \in M$ and $x \in \overline{\mathcal{P}(p)}$, and consequently
$$T_x M = \overline{\text{Lie}_x \mathcal{F}} = \overline{\mathcal{T}_x \mathcal{P}(p)}.$$
The theorem then follows from Corollary 3.2 with $\mathcal{B} = \overline{\mathcal{P}(p)}$. Q.E.D.

Here is an example to show how Theorem 3.1 works on the group of orientation-preserving diffeomorphisms of the unit circle.

**Example 3.3.** Let $S^1$ be the unit circle embedded into the Euclidean space $\mathfrak{R}^2$, and $M = \text{Diff}_0(S^1)$ be the identity connected component of the group of diffeomorphisms of $S^1$. In fact $M$ is a Lie group modelled on the locally convex space $\text{Vec}(S^1)$, cf. (Milnor 1984, pp 1039–41). Hence the tangent space of $M$ at $\text{id} \in M$ can be identified with
$$T_{\text{id}} M = \text{Vec}(S^1) = \{v(\theta)\partial_\theta \mid \theta \in S^1 (\text{mod } 2\pi), v \in C^\infty(S^1, \mathfrak{R})\},$$
where $\partial_\theta$ stands for $\frac{\partial}{\partial \theta}$. Under this identification, the commutator of two elements in the Lie algebra $\text{Vec}(S^1)$ of smooth vector fields on the circle is given by
$$[v(\theta)\partial_\theta, \omega(\theta)\partial_\theta] = (v'(\theta)\omega(\theta) - \omega'(\theta)v(\theta))\partial_\theta,$$
where $v'$ denotes the $\theta$-derivative of $v$. Note that this Lie bracket is the negative of the commonly assumed commutator of vector fields. Let $\text{Vec}(S^1)_\mathbb{C} = \text{Vec}(S^1) \otimes \mathbb{C}$ be the complexification of the Lie algebra $\text{Vec}(S^1)$. An element $v(\theta)\partial_\theta \in \text{Vec}(S^1)_\mathbb{C}$ can in fact be expressed using the Fourier expansion of $v(\theta) = \sum_{n=-\infty}^{+\infty} a_n e^{in\theta}$, where $a_n \in \mathfrak{R}$ and $e^{in\theta} = \cos n\theta + i \sin n\theta$. Hence $B_{\text{id}} := \{\partial_\theta, \cos n\theta \partial_\theta, \sin n\theta \partial_\theta\}_{n=1}^\infty$ forms a basis for $T_{\text{id}} M = \text{Vec}(S^1)$. Let $\tilde{B}_{\text{id}} = \{\cos\theta\partial_\theta, \sin\theta\partial_\theta, \cos 2\theta\partial_\theta, \sin 2\theta\partial_\theta\} \subset B_{\text{id}}$. It is easily seen that
$$[\sin\theta\partial_\theta, \cos\theta\partial_\theta] = \partial_\theta,$$
$$[ie^{in\theta}\partial_\theta, ie^{im\theta}\partial_\theta] = (m-n)ie^{i(m+n)\theta}\partial_\theta.$$
Comparing the real and imaginary parts of both sides of the latter equality, we deduce that taking linear combinations of all possible (iterated) Lie brackets of elements in $\tilde{B}_{\text{id}}$ one can generate all vector fields in $B_{\text{id}}$; e.g.,
$$\sin 3\theta\partial_\theta = [\cos\theta\partial_\theta, \cos 2\theta\partial_\theta] - [\sin\theta\partial_\theta, \sin 2\theta\partial_\theta],$$
$$\cos 3\theta\partial_\theta = -([\sin\theta\partial_\theta, \cos 2\theta\partial_\theta] + [\cos\theta\partial_\theta, \sin 2\theta\partial_\theta]),$$
$$\sin 4\theta\partial_\theta = ([\cos\theta\partial_\theta, \cos 3\theta\partial_\theta] - [\sin\theta\partial_\theta, \sin 3\theta\partial_\theta])/2, \text{ etc.}$$
Let us now consider $\tilde{B}_\phi := d_{\text{id}} R_\phi(\tilde{B}_{\text{id}}) \subset T_\phi M$, where $\phi \in M = \text{Diff}_0(S^1)$, and $R_\phi : M \ni \psi \mapsto \psi \circ \phi \in M$ is the right translation map. Accordingly, we can define the distribution $\mathcal{H} = \coprod_{\phi \in M} \mathcal{H}_\phi \subset TM$ where $\mathcal{H}_\phi := \text{span}\,\tilde{B}_\phi \subset T_\phi M$. Setting $\mathcal{F} := \{X \in C^{\infty}(M, \mathcal{H}) \mid X(\phi) \in \mathcal{H}_\phi \text{ for any } \phi \in M\}$, we conclude that

$$Lie_\phi \mathcal{F} = T_\phi M, \text{ for any } \phi \in M.$$

Theorem 3.1 now shows that $\overline{\mathcal{P}(\mathcal{F})(\phi)} = M$ for any $\phi \in M$.

It is worth noting that our controllability result on the group of diffeomorphisms $\text{Diff}_0(S^1)$, given in Example 3.3, does not follow from those obtained by Agrachev and Caponigro (2009).

## Acknowledgments

The authors gratefully acknowledge partial support by the grants NFR-204726/V30 and NFR-213440/BG from the Research Council of Norway during the preparation of this paper at the Department of Mathematics, University of Bergen.

E−mail : salehani.k h.m@gmail.com (M.K. Salehani)